\documentclass[english]{article}
\usepackage[margin=2.5cm]{geometry}
\usepackage[hyphens]{url}
\usepackage{hyperref}
\usepackage[hyphenbreaks]{breakurl}
\usepackage{amsfonts,amssymb,amsmath,amsthm,color}
\usepackage{graphicx}
\usepackage{bbm}
\usepackage{tikz}

\newtheorem{theorem}{Theorem}[section]
\newtheorem{proposition}[theorem]{Proposition}
\newtheorem{lemma}[theorem]{Lemma}

\newtheorem{fact}[theorem]{Fact}

\newtheorem{corollary}[theorem]{Corollary}
\newtheorem{definition}[theorem]{Definition}

\theoremstyle{definition}

\numberwithin{equation}{section}
\numberwithin{figure}{section}
\allowdisplaybreaks

\usepackage{babel}
\usepackage[T1]{fontenc}
\usepackage[utf8]{inputenc}

\newcommand{\eps}{\varepsilon}
\newcommand{\su}{\subseteq}

\newcommand{\Fp}{\mathbb{F}_p}

\newcommand{\Z}{\mathbb{Z}}

\usepackage{xcolor}

\title{Improving Behrend's construction:\\ Sets without arithmetic progressions in integers\\
and over finite fields}

\author{Christian Elsholtz\thanks{Institute for Analysis and Number Theory, TU Graz, Austria. Email: \url{elsholtz@math.tugraz.at}.
Supported by the joint FWF-ANR project Arithrand: FWF: I 4945-N and ANR-20-CE91-0006, and FWF grant DOC 183.}\and
Zach Hunter\thanks{Department of Mathematics, ETH Z\"urich, Switzerland. Email: \url{zach.hunter@math.ethz.ch}.}\and
Laura Proske\thanks{TU Graz, Austria. Email: \url{proske@student.tugraz.at}.}\and
Lisa Sauermann\thanks{Institute for Applied Mathematics, University of Bonn, Germany. Email: \url{sauermann@iam.uni-bonn.de}. Supported by the DFG Heisenberg Program.}}
\date{}

\begin{document}

\maketitle

\begin{abstract}\noindent
We prove new lower bounds on the maximum size of subsets $A\subseteq \{1,\dots,N\}$ or $A\subseteq \mathbb{F}_p^n$ not containing three-term arithmetic progressions. In the setting of $\{1,\dots,N\}$, this is the first improvement upon a classical construction of Behrend from 1946 beyond lower-order factors (in particular, it is the first quasipolynomial improvement). In the setting of $\mathbb{F}_p^n$ for a fixed prime $p$ and large $n$, we prove a lower bound of $(cp)^n$ for some absolute constant $c>1/2$ (for $c = 1/2$, such a bound can be obtained via classical constructions from the 1940s, but improving upon this has been a well-known open problem).
\end{abstract}

\section{Introduction}

The questions of estimating the maximum possible sizes of subsets of $\{1,\dots,N\}$ and of $\mathbb{F}_p^n$ without three-term arithmetic progressions are among the most central problems in additive combinatorics. Let us denote the maximum possible size of such subsets by $r_3(N)$ and $r_3(\mathbb{F}_p^n)$, respectively. So, formally, $r_3(N)$ is the maximum possible size of a subset $A\su \{1,\dots,N\}$ such that there do not exist distinct $x,y,z\in A$ with $x+z=2y$, and similarly $r_3(\mathbb{F}_p^n)$ is the maximum possible size of a subset $A\su \mathbb{F}_p^n$ such that there do not exist distinct $x,y,z\in A$ with $x+z=2y$. The problem of estimating $r_3(N)$ was raised by Erd\H{o}s and Tur\'an \cite{erdos-turan} in 1936, and has been intensively studied since then. It also has connections to problems in communication complexity, see \cite{chandra-furst-lipton, hambardzumyan-et-al, lipmaa}. The problem for $\mathbb{F}_p^n$ has also been studied for several decades (see \cite{meshulam}, see also \cite[p.\ 142]{salem-spencer2} for a slightly different but related setting).

In a breakthrough result in 2017, Ellenberg and Gijswijt \cite{ellenberg-gijswijt} proved that for any prime $p\ge 3$, there is an upper bound of the form
\begin{equation}\label{eq-EG-bound}
r_3(\mathbb{F}_p^n)\le (c_pp)^n
\end{equation}
for some constant $c_p<1$ only depending on $p$ (and their constant $c_p$ converges to $0.841\dots$ for $p\to \infty$, see \cite[Eq.\ (4.11)]{blasiak-et-al}). In the integer setting, in a more recent breakthrough Kelley and Meka \cite{kelley-meka} proved the upper bound
\[r_3(N)\le N\cdot \exp(-c(\log N)^{1/12})\]
for all $N\ge 3$, for some absolute constant $c>0$. This drastically improved upon all the previous bounds, obtained over many decades by Roth \cite{roth}, Heath-Brown \cite{heath-brown}, Szemer\'{e}di \cite{szemeredi-1990}, Bourgain \cite{bourgain-1999,bourgain-2008}, Sanders \cite{sanders2012,sanders2011} and Bloom--Sisask \cite{bloom-sisask-logarithmic}. Afterwards, using a modification of their method, this bound was improved to
\[r_3(N)\le N\cdot \exp(-c(\log N)^{1/9})\]
by Bloom and Sisask \cite{bloom-sisask}.

These upper bounds for $r_3(N)$ match the shape of a classical lower bound for this problem due to Behrend \cite{behrend} from 1946, which is of the form
\begin{equation}\label{eq-behrend}
r_3(N)\ge N\cdot 2^{-(2\sqrt{2}+o(1)) \sqrt{\log_2 N}}.
\end{equation}
Over the past almost eighty years, only the $o(1)$-term in this bound has been improved. In Behrend's original bound, this $o(1)$-term in the exponent encapsulated a factor of $(\log_2 N)^{-1/4}$, so the bound was of the form $r_3(N)\ge \Omega(N\cdot 2^{-2\sqrt{2}\sqrt{\log_2 N}}\cdot(\log_2 N)^{-1/4})$. In 2010, Elkin \cite{elkin} improved this factor to $(\log_2 N)^{1/4}$ instead, and an alternative proof for the bound with this improved $o(1)$-term was found by Green and Wolf \cite{green-wolf}.

In this paper, we give the first improvement to Behrend's \cite{behrend} classical lower bound beyond the $o(1)$-term in (\ref{eq-behrend}). As stated in the following theorem, we show that the constant factor $2\sqrt{2}\approx 2.828$ in the exponent can be improved to $2\sqrt{\log_2(24/7)}\approx 2.667$, proving that the classical bound (\ref{eq-behrend}) is not tight.

\begin{theorem}\label{thm-integers}
We have
\[r_3(N)\ge N\cdot 2^{-(C+o(1)) \sqrt{\log_2 N}}\]
with $C=2\sqrt{\log_2(24/7)}<2\sqrt{2}$.
\end{theorem}

Our proof of Theorem \ref{thm-integers} is motivated by studying three-term progression free sets in $\mathbb{F}_p^n$, for a fixed relatively large prime $p$ and large $n$. In this setting, one can adapt Behrend's construction \cite{behrend} (as noted by Tao and Vu in their book on additive combinatorics \cite[Exercise 10.1.3]{tao-vu} and also observed by Alon, see \cite[Lemma 17]{fox-pham}) to show
\begin{equation}\label{eq-Fpn-previous-bound}
r_3(\mathbb{F}_p^n)\ge \Big(\frac{p+1}{2}\Big)^{n-o(n)}
\end{equation}
for any fixed prime $p$ and large $n$. Alternatively, such a bound can also be shown via an adaptation of an earlier construction in the integer setting due to Salem and Spencer \cite{salem-spencer} from 1942 (see \cite[Theorem 2.13]{alon-shpilka-umans}). The asymptotic notation $o(n)$ in the bound (\ref{eq-Fpn-previous-bound}) is for $n\to \infty$ with $p$ fixed. The best quantitative bound for the $o(n)$-term in this statement is due to relatively recent work of the first author and Pach \cite[Theorem 3.10]{elsholtz-pach}, but beyond the $o(n)$-term, this bound has not been improved (except for specific small primes $p$, see the discussion below).

Comparing the upper and lower bounds for $r_3(\mathbb{F}_p^n)$ for a fixed (reasonably large) prime $p$ and large $n$ in (\ref{eq-EG-bound}) and (\ref{eq-Fpn-previous-bound}), there is still a large gap. Both of these bounds are roughly of the form $(cp)^n$ with $0<c<1$, but with different values of $c$. For the upper bound, the best known constant due to Ellenberg--Gijswijt \cite{ellenberg-gijswijt} is $c\approx 0.85$ (when the fixed prime is large enough). For the lower bound the constant $c=1/2$ from Behrend's construction \cite{behrend} or alternatively the Salem--Spencer construction \cite{salem-spencer} has not been improved in more than eighty years despite a lot of attention, especially after the upper bound of Ellenberg--Gijswijt appeared (see e.g. the blog post \cite{ellenberg-blog} as well as the discussion thereafter). Here, we finally improve this constant in the lower bound to be strictly larger than $1/2$.

\begin{theorem}\label{thm-finite-fields}
There is a constant $c>1/2$ such that for every prime $p$ and every sufficiently large positive integer $n$ (sufficiently large in terms of $p$), we have $r_3(\mathbb{F}_p^n)\ge (cp)^n$.
\end{theorem}

Breaking the barrier of $1/2$ in this result for $r_3(\mathbb{F}_p^n)$ relies on the same key insights as our lower bounds for $r_3(N)$ in Theorem \ref{thm-integers} improving Behrend's construction. The problem for $\mathbb{F}_p^n$ (the so-called finite field model) was popularized in the hope that improved upper bounds for $r_3(\mathbb{F}_p^n)$ would also lead to better upper bounds for $r_3(N)$ (see e.g. \cite{green-survey,peluse-survey,wolf-survey}). It turned out, however, that for $\mathbb{F}_p^n$ much stronger bounds hold than in the integer setting (for fixed $p$, the Ellenberg--Gijswijt bound stated in (\ref{eq-EG-bound}) is of the form $r_3(\mathbb{F}_p^n)\le (p^n)^{1-\gamma_p}$ for some $\gamma_p>0$, but in the integer setting one cannot hope for a bound of the form $r_3(N)\le N^{1-\gamma}$ for some fixed $\gamma>0$, as shown by Salem--Spencer \cite{salem-spencer} and Behrend \cite{behrend}). Nevertheless, this paper establishes a connection between lower bounds for $r_3(N)$ and for $r_3(\mathbb{F}_p^n)$. 

Our proof of Theorem \ref{thm-finite-fields} shows that one can take any $c< \sqrt{7/24}$, for example $c=0.54$. Even though this may not seem like a large improvement over $1/2$, it is the first qualitative improvement over the constant $1/2$ from the constructions of Salem--Spencer and Behrend from the 1940's. Both of these constructions lead to three-term progression free subsets of $\mathbb{F}_p^n$ only consisting of vectors with all entries in $\{0,1,\dots,(p-1)/2\}$, i.e.\ they only use roughly half of the available elements in $\Fp$ in each coordinate. The restriction of all entries to $\{0,1,\dots,(p-1)/2\}$ is crucial in these constructions, as it ensures that there is no ``wrap-around'' over $\Fp$. However, such an approach cannot be used to obtain three-term progression free subset $A\su \mathbb{F}_p^n$ of size $|A|> ((p+1)/2)^n$, so $c=1/2$ is a significant barrier for this problem. In light of this, it may actually be considered a surprise that it is possible to obtain a constant $c>1/2$.

For very small primes, better bounds were obtained with specific constructions depending on the particular prime. For example, in $\mathbb{F}_3^n$ a lower bound of $2.2202^{n-o(n)}$ was recently obtained by Romera-Paredes et al.\ \cite{google-deep-mind} using artifical intelligence building upon traditional methods from previous bound \cite{calderbank-fishburn,edel,tyrrell}. Naslund \cite{naslund} informed us about forthcoming work with an approach related to Shannon capacity, proving a lower bound of $2.2208^{n-o(n)}$ for the maximum size of three-term progression free subsets of $\mathbb{F}_3^n$. In $\mathbb{F}_5^n$, the best known lower bounds is $(35^{1/3})^{n-o(n)}$ due to the first and third author, Pollak, Lipnik and Siebenhofer (see \cite{elsholtz-lipnik-siebenhofer}), note that $35^{1/3}\approx 3.271$.

Our methods for improving the lower bounds for $r_3(N)$ and $r_3(\mathbb{F}_p^n)$ can be applied more generally for finding three-term progression free subsets in finite abelian groups. We recall that a three-term arithmetic progression in a finite abelian group $G$ consists of three distinct elements $x,y,z\in G$ with $x+z=y+y$ (see e.g. \cite{croot-lev-pach}), and let us denote the maximum possible size of a three-term progression free subset of $G$ by $r_3(G)$. Brown and Buhler \cite{brown-buhler} proved in 1982 that $r_3(G) \leq o(|G|)$ for all finite abelian groups $G$ (where the asymptotic notation $o(|G|)$ is for growing group size $|G|$). Frankl, Graham and R\"odl \cite{frankl-graham-rodl} gave a different proof of this statement in 1987, and further upper bounds were obtained by Meshulam \cite{meshulam} and Lev \cite{lev:2004}.

Of course, every finite abelian group $G$ can be written in the form $\Z_{m_1}\times \dots\times \Z_{m_n}$, and the problem of estimating $r_3(\Z_{m_1}\times \dots\times \Z_{m_n})$ has been most intensively studied in the case $m_1=\dotsm=m_n$, i.e. in the setting of $\Z_m^n$. In breakthrough work, Croot--Lev--Pach \cite{croot-lev-pach} proved the upper bound $r_3(\Z_4^n)\le 3.611^n$ by introducing a new polynomial method, which also led to the Ellenberg--Gijswijt bound for $r_3(\mathbb{F}_p^n)$ stated in (\ref{eq-EG-bound}). Using similar methods, Petrov and Pohoata \cite{petrov-pohoata} studied upper bounds for $r_3(\Z_8^n)$.

Generalizing Theorem \ref{thm-finite-fields}, we show that for any integer $m\ge 2$ (not necessarily prime) and sufficiently large $n$, we have $r_3(\Z_m^n)\ge (cm)^n$ for some absolute constant $c>1/2$. For odd $m$, the best previous lower bound for this problem has (as in the case of prime $m$) been of the form $((m+1)/2)^{n-o(n)}$ (based on the Salem--Spencer construction \cite{salem-spencer} or the Behrend construction \cite{behrend}). For even $m$, the best previous lower bound has been $((m+2)/2)^{n-o(n)}$ due to the first author and Pach \cite[Theorem 3.11]{elsholtz-pach}.

In fact, we show the following more general lower bound for arbitrary abelian groups $\Z_{m_1}\times \dots\times \Z_{m_n}$.

\begin{theorem}\label{thm-Zm}
Consider integers $m_1,\dots,m_n$ (for some positive integer $n$), and let $m$ be such that $2\le m_i\le m$ for $i=1,\dots,n$. Then we have
\[r_3(\Z_{m_1}\times \dots\times \Z_{m_n})\ge \frac{(7/24)^{n/2}}{10^{6}m^2n^3}\cdot m_1\dotsm m_n.\]
\end{theorem}

Note that Theorem \ref{thm-Zm} immediately implies that for any fixed constant $c<\sqrt{7/24}$ and for every integer $m\ge 2$, we have $r_3(\Z_m^n)\ge (cm)^n$ if $n$ is sufficiently large with respect to $m$ and $c$. In particular, it implies Theorem \ref{thm-finite-fields}.

It is plausible that by slight modifications of our method, one can obtain better numerical bounds in our results above. In particular, our construction relies on certain explicit two-dimensional building blocks of area close to $7/24$, and an improved choice for these building blocks with larger area would automatically carry over to numerical improvements of the constant $C=2\sqrt{\log_2(24/7)}$ in Theorem \ref{thm-integers} and the constant $7/24$ in Theorem \ref{thm-Zm}. We see the main contribution of this paper as introducing this method, and using it to break the lower bounds for $r_3(N)$ and $r_3(\mathbb{F}_p^n)$ in (\ref{eq-behrend}) and (\ref{eq-Fpn-previous-bound}) originating from the 1940's.

In the next section we give an overview of our proof approach and state some key propositions, from which we will deduce Theorems \ref{thm-integers} to \ref{thm-Zm} in Section \ref{sect-deduction-main-results}.  To prove these propositions, we need to find suitable two-dimensional building blocks, and these will be constructed in Section \ref{sect-constr-building-block}.

\textit{Notation.} As usual, for a point $x\in \mathbb{R}^n$, we denote the coordinates of $x$ by $x_1,\dots,x_n$. For real numbers $a,a'\in \mathbb{R}$, we write  $a\equiv a'\mod 1$ if $a-a'\in \mathbb{Z}$. For two points $x\in \mathbb{R}^n$ and $y\in \mathbb{R}^n$, we write $x\equiv y\mod 1$ if we have $x_i\equiv y_i\mod 1$ for all $i=1,\dots,n$ (i.e.\ if we have $x-y\in \mathbb{Z}^n$). For a measurable subset $S\su \mathbb{R}^n$, we denote its measure by $\mu(S)$.

\textit{Acknowledgements.} We would like to thank Daniel Carter, Zachary Chase, Ben Green, Eric Naslund, and Benny Sudakov for interesting discussions and helpful comments.

\section{Proof Overview}
\label{sect-overview}

The first step in our constructions for large three-term progression free subsets of $\{1,\dots,N\}$ or $\Z_{m_1}\times \dots\times \Z_{m_n}$ is to take a (randomized) embedding into a (high-dimensional) torus. This way, we can reduce both problems to finding a large subset of the torus without a three-term arithmetic progression such that the first and third term are reasonably far away from each other. Note that we cannot aim to find a large subset of the torus without any three-term arithmetic progressions, since it is well-known that any three-term progression free subset of the torus has measure zero (see for example \cite{ruziewicz}). So it is essential to only forbid three-term arithmetic progressions with far-away first and third term in our desired subset of the torus.

The following proposition states that there is indeed such a large-measure subset of a (high-dimensional) torus without three-term arithmetic progressions with far-away first and third term. Proving this proposition is the main difficulty of this paper.

\begin{proposition}\label{prop-main}
For any $0<\delta<1$ and any even positive integer $n$, there exists a measurable subset $S\su [0,1)^n$ with measure $\mu(S) \ge 10^{-5}\delta^2n^{-3}\cdot (7/24)^{n/2}$ such that for any $x,y,z\in S$ with $x+z\equiv 2y \mod 1$ we have $|x_i-z_i|<\delta$ for $i=1,\dots,n$.
\end{proposition}

The deductions of Theorems~\ref{thm-integers} to \ref{thm-Zm} from Proposition~\ref{prop-main} will be discussed in Section \ref{sect-deduction-main-results}.

We remark that it is easy to show a weaker version of Proposition \ref{prop-main} with a bound of $\mu(S) \ge \delta^2\cdot (1/2)^{n+o(n)}$ for the measure of $S$. Indeed, there one can take $S\su [0,1/2)^n$ to be an appropriately chosen sphere shell of width on the order of $\delta^2$ (notice that then for $x,y,z\in S$ we have $x+z\equiv 2y \mod 1$ if and only if  $x+z=2y$, since $S\su [0,1/2)^n$). From this weaker statement one can deduce (\ref{eq-behrend}) with Behrend's original constant $2\sqrt{2}$, this approach to proving Behrend's bound (with an improved $o(1)$-term) is due to Green--Wolf \cite{green-wolf}. Obtaining the stronger bound on $\mu(S)$ in Proposition~\ref{prop-main} (having a factor of $(7/24)^{n/2}$ instead of $(1/2)^{n}=(1/4)^{n/2}$), which leads to our improvement upon the constant $2\sqrt{2}$, is much more difficult. In particular, for $n$ large in terms of $\delta$, any subset $S\su [0,1)^n$ whose measure is as large as in Proposition \ref{prop-main} will have points $x,y,z\in S$ with $x+z\equiv 2y \mod 1$ but $x+z\ne 2y$, so the ``modulo 1'' cannot just be ignored like in the previous approach. In fact, there will always be $x,y,z\in S$ with $x+z\equiv 2y \mod 1$ and $x=z\ne y$. Therefore, in the conclusion of Proposition \ref{prop-main}, it is crucial to bound the absolute values of the coordinates of $x-z$ and not of $x-y$ (which, as the step length of the ``modulo 1'' arithmetic progression $x,y,z$, would maybe be more natural to consider), since the analogous statement for $x-y$ would be false.

In order to prove Proposition~\ref{prop-main}, we use subsets of a two-dimensional torus as building blocks. These two-dimensional building blocks need to satisfy a somewhat technical condition, summarized in the following proposition.

\begin{proposition}\label{prop-building-block}
Let $0<\eps<1$. Then there exists a measurable subset $T\su [0,1)^2$ with measure $\mu(T) \ge 7/24-\eps$ and a measurable function $f: T\to [0,100/\eps^2]$, such that for any $x,y,z\in T$ with $x+z\equiv 2y \mod 1$ we have
\[f(x) + f(z) \ge 2f(y) + (x_1-z_1)^2 + (x_2-z_2)^2.\]
\end{proposition}

To deduce Proposition~\ref{prop-main} from Proposition~\ref{prop-building-block}, we define the set $S\su [0,1)^n$ in Proposition~\ref{prop-main} to be a subset of the $(n/2)$-fold product set $T^{n/2}\su [0,1)^n$ consisting of those points $((x_1,x_2),(x_3,x_4),\dots,(x_{n-1},x_n))\in T^{n/2}$ where $f(x_1,x_2)+f(x_3,x_4)+\dots+f(x_{n-1},x_n)$ is contained in some small interval. In other words, the set $S\su [0,1)^n$ can be viewed as a ``slice'' of the product set $T^{n/2}\su [0,1)^n$ with respect to the function $f(x_1,x_2)+f(x_3,x_4)+\dots+f(x_{n-1},x_n)$. It is then not hard to show, using the conditions on the set $T$ and the function $f$ in Proposition~\ref{prop-building-block}, that such a set $S$ satisfies the conditions in Proposition~\ref{prop-main}. The details of this proof are given at the end of this section.

In order to prove Proposition~\ref{prop-building-block}, we explicitly construct the desired subset $T\su [0,1)^2$ and the function $f$. The set $T$ is a carefully chosen union of polygons in $[0,1)^2$. The construction is somewhat ad-hoc and it seems plausible that the constant $7/24$ in Propositions~\ref{prop-building-block} and \ref{prop-main} (and consequently also the constants in Theorems~\ref{thm-integers} and \ref{thm-Zm}) can be improved by finding a better construction of such a set $T$. However, the main novelty of this paper is our overall approach for proving Theorems~\ref{thm-integers} to \ref{thm-Zm}, and showing that this approach can be used to beat the long-standing bounds on these problems.

The proof of Proposition~\ref{prop-building-block} can be found in Section \ref{sect-constr-building-block}. We finish this section by showing that Proposition~\ref{prop-building-block} indeed implies Proposition~\ref{prop-main}.

\begin{proof}[Proof of Proposition~\ref{prop-main} assuming Proposition~\ref{prop-building-block}] For $n=2$, Proposition~\ref{prop-main} is trivially true (for example, we can take $S=[0,\delta)^2$). So we may assume that $n\ge 4$.

Define $\eps=1/n$ and let $T\su [0,1)^2$ and the function $f$ be as in Proposition~\ref{prop-building-block}. Note that then the $(n/2)$-fold product set $T^{n/2}\su [0,1)^n$ has measure
\[\mu(T^{n/2})=\mu(T)^{n/2}\ge (7/24-\eps)^{n/2}\ge \Big(1-\frac{24}{7}\eps\Big)^{n/2}\cdot \Big(\frac{7}{24}\Big)^{n/2}\ge  10^{-2}\cdot \Big(\frac{7}{24}\Big)^{n/2},\]
where we used that for all $a\in [0,6/7]$ we have $1-a\ge 10^{-(7/6)a}$ (indeed, for $a=0$ we have $1-0= 10^{-(7/6)\cdot 0}$, for $a=6/7$ we have $1-(6/7)=1/7>1/10= 10^{-(7/6)\cdot (6/7)}$, and the function $10^{-(7/6)a}$ is concave on the interval $[0,6/7]$) and therefore (recalling that $
\eps=1/n\le 1/4$)
\[\Big(1-\frac{24}{7}\eps\Big)^{n/2}\ge \Big(10^{-(7/6)\cdot (24/7)\eps}\Big)^{n/2}=10^{-2\eps n}=10^{-2}.\]
On the set $T^{n/2}$, let us consider the function taking value $f(x_1,x_2)+f(x_3,x_4)+\dots+f(x_{n-1},x_n)$ for each $((x_1,x_2),(x_3,x_4),\dots,(x_{n-1},x_n))\in T^{n/2}$. This is a measurable function, with values in the interval $[0,100n^3]$ (noting that $n\cdot 100/\eps^2=100n^3$). Thus, for every integer $j=0,\dots,\lfloor 100n^3\cdot (2/\delta^2)\rfloor $, the set 
\[S_j=\{((x_1,x_2),\dots,(x_{n-1},x_n))\in T^{n/2}\mid j\cdot \delta^2/2\le f(x_1,x_2)+\dots+f(x_{n-1},x_n)<(j+1)\cdot \delta^2/2\}\]
is measurable, and the union of these sets $S_j$ for $j=0,\dots,\lfloor 100n^3\cdot (2/\delta^2)\rfloor $ is the entire product set $T^{n/2}$. Thus, there exists some $j\in \{0,\dots,\lfloor 100n^3\cdot (2/\delta^2)\rfloor\}$ such that
\[\mu(S_j)\ge \frac{\mu(T^{n/2})}{\lfloor 100n^3\cdot (2/\delta^2)\rfloor+1}\ge \frac{10^{-2}\cdot (7/24)^{n/2}}{4\cdot 100n^3\delta^{-2}}\ge 10^{-5}\delta^2n^{-3}\cdot (7/24)^{n/2}.\]
It remains to show that for any $x,y,z\in S_j$ with $x+z\equiv 2y \mod 1$ we have $|x_i-z_i|<\delta$ for $i=1,\dots,n$. Indeed, for every $h=1,\dots,n/2$ we have $(x_{2h-1},x_{2h}),(y_{2h-1},y_{2h}),(z_{2h-1},z_{2h})\in T$ and $(x_{2h-1},x_{2h})+(z_{2h-1},z_{2h})\equiv 2(y_{2h-1},y_{2h})\mod 1$ and so by the condition in Proposition~\ref{prop-building-block} we obtain
\[f(x_{2h-1},x_{2h}) + f(z_{2h-1},z_{2h}) \ge 2f(y_{2h-1},y_{2h}) + (x_{2h-1}-z_{2h-1})^2 + (x_{2h}-z_{2ih})^2.\]
Summing this up for $h=1,\dots,n/2$, and recalling the definition of $S_j$, yields
\begin{align*}
2\cdot (j+1)\cdot \delta^2/2&> f(x_1,x_2)+\dots+f(x_{n-1},x_n)+f(z_1,z_2)+\dots+f(z_{n-1},z_n)\\
&\ge 2 f(y_1,y_2)+\dots+2f(y_{n-1},y_n)+ (x_{1}-z_{1})^2 + \dots+ (x_{n}-z_{n})^2\\
&\ge 2\cdot j\cdot \delta^2/2+(x_{1}-z_{1})^2 + \dots+ (x_{n}-z_{n})^2,
\end{align*}
and hence
\[(x_{1}-z_{1})^2 + \dots+ (x_{n}-z_{n})^2< 2\cdot \delta^2/2=\delta^2.\]
This implies that $|x_i-z_i|<\delta$ for $i=1,\dots,n$, as desired.
\end{proof}

\section{Deduction of main results}
\label{sect-deduction-main-results}

In this section we derive Theorems~\ref{thm-integers} to \ref{thm-Zm} from Proposition~\ref{prop-main}. We start with Theorem \ref{thm-Zm}, and first prove the following version of it for the case of even $n$.

\begin{proposition}\label{prop-Zm-even-n}
Consider integers $m_1,\dots,m_n$ (for some positive integer $n$), and let $m$ be such that $2\le m_i\le m$ for $i=1,\dots,n$. If $n$ is even, we have
\[r_3(\Z_{m_1}\times \dots\times \Z_{m_n})\ge \frac{(7/24)^{n/2}}{10^{5}m^2n^3}\cdot m_1\dotsm m_n.\]
\end{proposition}

As mentioned in the previous section, this statement is deduced from Proposition~\ref{prop-main} by considering a (randomized) embedding of $\Z_{m_1}\times \dots\times \Z_{m_n}$ into the torus $[0,1)^n$. For a subset $S\su [0,1)^n$ as in Proposition~\ref{prop-main}, one can then show that its pre-image is a three-term progression free subset. Using the lower bound on $\mu(S)$ in Proposition~\ref{prop-main}, this gives a lower bound for $r_3(\Z_{m_1}\times \dots\times \Z_{m_n})$.

\begin{proof}
Let $\delta=1/m$, and let $S\su [0,1)^n$ be a measurable subset as in Proposition~\ref{prop-main}. For any $a\in [0,1]^n$, let us now define a map $\varphi_{a}:\Z_{m_1}\times \dots\times \Z_{m_n}\to [0,1)^n$ as follows. Every point in $\Z_{m_1}\times \dots\times \Z_{m_n}$ can be represented (uniquely) by an $n$-tuple $(r_1,\dots,r_n)\in \{0,1,\dots,m_1-1\}\times \dots \times \{0,1,\dots,m_n-1\}$. Let us now define $\varphi_{a}(r_1,\dots,r_n)=(q_1,\dots,q_n)$ for the unique point $q=(q_1,\dots,q_n)\in [0,1)^n$ with $(q_1,\dots,q_n)\equiv (a_1+r_1/m_1,\dots, a_n+r_n/m_n)\mod 1$.

Note that the map $\varphi_{a}:\Z_{m_1}\times \dots\times \Z_{m_n}\to [0,1)^n$ is injective for any $a\in [0,1]^n$. Furthermore, note that for any $a\in [0,1]^n$ and any $x,z\in \Z_{m_1}\times \dots\times \Z_{m_n}$ we either have $\varphi_{a}(x)_i=\varphi_{a}(z)_i$ or $|\varphi_{a}(x)_i-\varphi_{a}(z)_i|\ge 1/m_i\ge 1/m=\delta$ (indeed, the different values occurring as the $i$-th coordinate of points in the image $\varphi_{a}(\Z_{m_1}\times \dots\times \Z_{m_n})$ form an arithmetic progression with step-width $1/m_i$). Hence for any $a\in [0,1]^n$ and any distinct $x,z\in \Z_{m_1}\times \dots\times \Z_{m_n}$, there exists a coordinate $i\in \{1,\dots,n\}$ such that $|\varphi_{a}(x)_i-\varphi_{a}(z)_i|\ge \delta$.

We claim that for any $a\in [0,1]^n$, the pre-image $\varphi_{a}^{-1}(S)$ of the set $S$ from Proposition \ref{prop-main} is a three-term progression free subset of $\Z_{m_1}\times \dots\times \Z_{m_n}$. Indeed, suppose that for some $a\in [0,1]^n$, there exist distinct $x,y,z\in \varphi_{a}^{-1}(S)\su \Z_{m_1}\times \dots\times \Z_{m_n}$ with $x+z=2y$. Then, as shown above, there must be a coordinate $i\in \{1,\dots,n\}$ such that  $|\varphi_{a}(x)_i-\varphi_{a}(z)_i|\ge \delta$. On the other hand, having $x+z=2y$ means that $x_i+z_i\equiv 2y_i\mod m_i$ for all $i=1,\dots,n$, where $(x_1,\dots,x_n),(y_1,\dots,y_n), (z_1,\dots,z_n)\in \{0,1,\dots,m_1-1\}\times \dots \times \{0,1,\dots,m_n-1\}$ are the $n$-tuples corresponding to $x,y,z\in \Z_{m_1}\times \dots\times \Z_{m_n}$, respectively. Hence we obtain $(x_i/m_i)+(z_i/m_i)\equiv 2(y_i/m_i)\mod 1$ for $i=1,\dots,n$, and consequently
\[\varphi_{a}(x)_i+\varphi_{a}(z)_i\equiv (a_i+x_i/m_i)+(a_i+z_i/m_i)\equiv 2(a_i+y_i/m_i)\equiv 2\varphi_{a}(y)_i\mod 1\]
for $i=1,\dots,n$. This means that $\varphi_{a}(x)+\varphi_{a}(z)\equiv 2\varphi_{a}(y) \mod 1$. Since $\varphi_{a}(x),\varphi_{a}(y),\varphi_{a}(z)\in S$, the condition for the set $S$ in Proposition \ref{prop-main} now implies that $|\varphi_{a}(x)_i-\varphi_{a}(z)_i|< \delta$ for all $i=1,\dots,n$. This is a contradiction, so the pre-image $\varphi_{a}^{-1}(S)\su \Z_{m_1}\times \dots\times \Z_{m_n}$ is indeed three-term progression free for any $a\in [0,1]^n$.

Finally, consider a uniformly random choice of $a\in [0,1]^n$. For any $x\in \Z_{m_1}\times \dots\times \Z_{m_n}$ we have $x\in \varphi_{a}^{-1}(S)$ if and only if $\varphi_{a}(x)\in S$, and this happens with probability $\mu(S)$. Hence
\[\mathbb{E}[|\varphi_{a}^{-1}(S)|]=\sum_{x\in \Z_{m_1}\times \dots\times \Z_{m_n}}\mathbb{P}[x\in \varphi_{a}^{-1}(S)]=|\Z_{m_1}\times \dots\times \Z_{m_n}|\cdot \mu(S)\ge \frac{(7/24)^{n/2}}{10^{5}m^2n^3}\cdot m_1\dotsm m_n,\]
recalling that $\mu(S) \ge 10^{-5}\delta^2n^{-3}\cdot (7/24)^{n/2}=10^{-5}(1/m)^2n^{-3}\cdot (7/24)^{n/2}$. This implies that for some choice of $a\in [0,1]^n$ we must have
\[|\varphi_{a}^{-1}(S)|\ge \frac{(7/24)^{n/2}}{10^{5}m^2n^3}\cdot m_1\dotsm m_n.\]
As shown above, the set $\varphi_{a}^{-1}(S)\su \Z_{m_1}\times \dots\times \Z_{m_n}$ is three-term progression free, so this establishes the desired lower bound for $r_3(\Z_{m_1}\times \dots\times \Z_{m_n})$.
\end{proof}

To deduce Theorem \ref{thm-Zm} from Proposition \ref{prop-Zm-even-n}, we use the following simple fact.

\begin{fact}\label{fact-product}
    For any finite abelian groups $G$ and $H$, we have
    \[r_3(G)\ge \frac{r_3(G \times H)}{|H|}.\]
\end{fact}

\begin{proof}
    Let $A\su G\times H$ be a three-term progression free subset of $G\times H$ of size $|A|=r_3(G\times H)$. For any $h\in H$, let us define the map $\varphi_h:G\to  G\times H$ by $\varphi_h(g)=(g,h)$ for all $g\in G$. For any $h\in H$, this map is injective and for any distinct $x,y,z\in G$, the images $\varphi_h(x), \varphi_h(y), \varphi_h(z)\in G\times H$ form a three-term arithmetic progression in $G\times H$ if and only if $x,y,z$ form a three-term arithmetic progression in $G$. Thus, for any $h\in H$, the pre-image $\varphi_h^{-1}(A)\su G$ is a three-term progression free subset of $G$ and therefore has size $|\varphi_h^{-1}(A)|\le r_3(G)$. On the other hand, by injectivity of $\varphi_h$ we have $|\varphi_h^{-1}(A)|=|A\cap \varphi_h(G)|$. As the images $\varphi_h(G)$ for all $h\in H$ form a partition of $G\times H$, this implies
    \[r_3(G\times H)=|A|=\sum_{h\in H} |A\cap \varphi_h(G)|=\sum_{h\in H}|\varphi_h^{-1}(A)|\le |H|\cdot r_3(G).\]
    Rearranging yields the desired inequality.
\end{proof}

Now, we are ready to prove Theorem \ref{thm-Zm}

\begin{proof}[Proof of Theorem \ref{thm-Zm}]
Note that the inequality in the theorem statement holds trivially for $n\le 2$ (since then the right-hand side is smaller than $1$). So let us assume $n\ge 3$. If $n$ is even, then the desired lower bound for $r_3(\Z_{m_1}\times \dots\times \Z_{m_{n}})$ follows immediately from Proposition \ref{prop-Zm-even-n}. If $n$ is odd, applying Proposition \ref{prop-Zm-even-n} with $m_{n+1}=m$ yields
\[r_3(\Z_{m_1}\times \dots\times \Z_{m_n}\times \Z_{m})\ge \frac{(7/24)^{(n+1)/2}}{10^{5}m^2(n+1)^3}\cdot m_1\dotsm m_n\cdot m\ge \frac{(7/24)^{n/2}}{10^{6}m^2n^3}\cdot m_1\dotsm m_n\cdot m.\]
Applying Fact \ref{fact-product} to $G=\Z_{m_1}\times \dots\times \Z_{m_n}$ and $H=\Z_{m}$, we obtain
\[r_3(\Z_{m_1}\times \dots\times \Z_{m_{n}})\ge \frac{r_3(\Z_{m_1}\times \dots\times \Z_{m_{n}}\times \Z_{m})}{m}\ge \frac{(7/24)^{n/2}}{10^{6}m^2n^3}\cdot m_1\dotsm m_n.\qedhere\]
\end{proof}

Theorem \ref{thm-finite-fields} is an immediate consequence of Theorem \ref{thm-Zm}.

\begin{proof}[Proof of Theorem \ref{thm-finite-fields}]
Let $1/2<c<\sqrt{7/24}$. Then, as long as $n$ is sufficiently large with respect to $p$, by Theorem \ref{thm-Zm} we have
\[r_3(\mathbb{F}_p^n)=r_3(\Z_p\times\dots\times\Z_p)\ge \frac{(7/24)^{n/2}}{10^{6}p^2n^3}\cdot p^n=\frac{(\sqrt{7/24}\cdot p)^{n}}{10^{6}p^2n^3}\ge (cp)^n.\qedhere\]
\end{proof}

Finally, it remains to derive Theorem~\ref{thm-integers}.

\begin{proof}[Proof of Theorem~\ref{thm-integers}]
As we are proving an asymptotic statement, we may assume that $N$ is sufficiently large. Let $n$ be an even positive integer with
\begin{equation}\label{eq-choice-of-n}
\frac{2\cdot \sqrt{\log_2 N}}{\sqrt{\log_2(24/7)}}\le n\le \frac{2\cdot \sqrt{\log_2 N}}{\sqrt{\log_2(24/7)}}+2.
\end{equation}
Let $p_1,\dots, p_n$ be the first $n$ primes. Defining $p=\max\{p_1,\dots,p_n\}$, we have $p\le 100 n\log_2 n$ (this follows from the prime number theorem, but there are also easier proofs for this weaker bound, see e.g. \cite[p.\ 150]{sierpinski} for a relatively short elementary proof). In particular, we can observe that $p_1\dotsm p_n\le p^n\le (100n\log_2 n)^n\le (n^2/4)^{n}\le (\log_2 N)^{\sqrt{\log_2 N}}<N$.

We claim that there are powers $m_1,\dots,m_n$ of the primes $p_1,\dots,p_m$, respectively, with $N/p\le m_1\dotsm  m_n\le N$ and $p_i\le m_i< N^{1/n}\cdot p_i$ for $i=1,\dots,n$. Indeed, for each $i=1,\dots,n$, start by defining $m_i'$ to be the unique power of $p_i$ with $N^{1/n}\le m_i'<N^{1/n}\cdot p_i$. Note that then we have $m_1'\dotsm  m_n'\ge (N^{1/n})^n=N$, and $m_i'< N^{1/n}\cdot p$ for $i=1,\dots,n$. As long as $m_1'\dotsm  m_n'>N$, let us decrease the product $m_1'\dotsm  m_n'$ by choosing an index $i$ with $m_i'>p_i$ and dividing $m_i'$ by $p_i$ (again obtaining a power of $p_i$). In each step the product $m_1'\dotsm  m_n'$ decreases by a factor of at most $p$. Thus, when the product $m_1'\dotsm  m_n'$ stops being larger than $N$, it must attain a value between $N/p$ and $N$. We can then define $m_1,\dots,m_n$ to be the values of $m_1',\dots,m_n'$ at that point.

Note that $m_1,\dots,m_n$ are coprime (as they are powers of the distinct primes $p_1,\dots, p_n$), and hence $\mathbb{Z}_{m_1\dotsm m_n}\cong \Z_{m_1}\times \dots\times \Z_{m_{n}}$. Thus, for $m=\lfloor N^{1/n}\cdot p\rfloor$, we obtain
\[r_3(N)\ge r_3(m_1\dotsm  m_n)\ge r_3(\mathbb{Z}_{m_1\dotsm m_n})=r_3(\Z_{m_1}\times \dots\times \Z_{m_{n}})\ge \frac{(7/24)^{n/2}}{10^{5}m^2n^3}\cdot m_1\dotsm m_n,\]
where at the last step we used Proposition \ref{prop-Zm-even-n} (and at the first two steps we used that every three-term progression free subset of $\mathbb{Z}_{m_1\dotsm m_n}$ gives rise to a three-term progression free subset of $\{1,\dots,m_1\dots m_n\}$ and hence in particular to a three-term progression free subset of $\{1,\dots,N\}$, recalling that $m_1\dotsm  m_n\le N$). Rewriting the right-hand side yields
\[r_3(N)\ge \frac{(7/24)^{n/2}}{10^{5}m^2n^3}\cdot m_1\dotsm m_n\ge \frac{(7/24)^{n/2}}{10^{5}(N^{1/n}p)^2n^3}\cdot (N/p)=\frac{1}{10^{5}p^3n^3}\cdot N\cdot (7/24)^{n/2}\cdot N^{-2/n}.\]
Recalling $p\le 100 n\log_2 n$, and recalling our choice of $n$, we can conclude that
\begin{align*}
r_3(N)&\ge \frac{1}{10^{5}(100 n\log_2 n)^3n^3}\cdot N\cdot 2^{-\log_2(24/7)\cdot n/2}\cdot 2^{-2(\log_2 N)/n}\\
&\ge \frac{1}{10^{11}n^6(\log_2 n)^3}\cdot N\cdot 2^{-\sqrt{\log_2(24/7)}\cdot \sqrt{\log_2 N}-2}\cdot 2^{-\sqrt{\log_2(24/7)}\cdot \sqrt{\log_2 N}}\\
&\ge \frac{1}{10^{12}(2\sqrt{\log_2 N})^6(\log_2 \log_2  N)^3}\cdot N\cdot 2^{-2\sqrt{\log_2(24/7)}\cdot \sqrt{\log_2 N}}\\
&\ge \frac{N\cdot 2^{-2\sqrt{\log_2(24/7)}\cdot \sqrt{\log_2 N}}}{10^{14}(\log_2 N)^3(\log_2 \log_2  N)^3}. 
\end{align*}
This shows that $r_3(N)\ge N\cdot 2^{-(2\sqrt{\log_2(24/7)}+o(1)) \cdot \sqrt{\log_2 N}}$, as desired.
\end{proof}

We remark that it is also possible to deduce Theorem \ref{thm-integers} directly from Proposition \ref{prop-main}, eliminating the need to use estimates for the size of the $n$-th prime as in the proof above. Indeed, following an idea of Green and Wolf \cite{green-wolf}, given $N$, for any even positive integer $n$ and any $a,b\in [0,1]^n$, one can consider the map $\varphi_{a,b}:\{1,\dots,N\}\to [0,1)^n$ defined by $\varphi_{a,b}(x)\equiv a+xb\mod 1$ for all $x\in \{1,\dots,N\}$. For a uniformly random choice of $b\in [0,1]^n$, the probability of having some $t\in \{1,\dots,N\}$ such that $tb\equiv v\mod 1$ for a vector $v\in \mathbb{R}^n$ with $\Vert v\Vert_\infty\le N^{-1/n}/4$ is at most $N\cdot (N^{-1/n}/2)^n\le 1/2$ (indeed, for every $t\in \{1,\dots,N\}$ the vector $tb \mod 1$ is uniformly distributed over $\mathbb{R}^n/\mathbb{Z}^n$ and so the probability of having $tb\equiv v\mod 1$ for a vector $v\in \mathbb{R}^n$ with $\Vert v\Vert_\infty\le N^{-1/n}/4$ is $(N^{-1/n}/2)^n$). So we can fix some $b\in [0,1]^n$ such that $tb\not\equiv v\mod 1$ for any  $t\in \{1,\dots,N\}$ and any vector $v\in \mathbb{R}^n$ with $\Vert v\Vert_\infty\le N^{-1/n}/4$. Let us now apply Proposition \ref{prop-main} with $\delta=N^{-1/n}/4$, yielding a measurable subset $S\su [0,1)^n$ satisfying the conditions in the proposition. Now, for any $a\in [0,1]^n$ the pre-image $\varphi_{a,b}^{-1}(S)\su \{1,\dots,N\}$ is three-term progression free. Indeed, suppose there were $x,y,z\in \varphi_{a,b}^{-1}(S)$ with $x>y>z$ and $x+z= 2y$. Then we would have $\varphi_{a,b}(x)+\varphi_{a,b}(z)\equiv (a+xb)+ (a+zb)=2a+(x+z)b=2(a+yb)\equiv 2\varphi_{a,b}(y)\mod 1$ and hence $|\varphi_{a,b}(x)_i-\varphi_{a,b}(z)_i|<\delta$ for $i=1,\dots,n$ by the conditions on $S$ (noting that $\varphi_{a,b}(x), \varphi_{a,b}(y), \varphi_{a,b}(z)\in S$), so $\Vert\varphi_{a,b}(x)-\varphi_{a,b}(z)\Vert_\infty<\delta=N^{-1/n}/4$. On the other hand, $(x-z)b= (a+xb)-(a+zb)\equiv \varphi_{a,b}(x)-\varphi_{a,b}(z)\mod 1$, which yields a contradiction to our choice of $b$ (as $x-z\in \{1,\dots,N\}$). Thus, the pre-image $\varphi_{a,b}^{-1}(S)\su \{1,\dots,N\}$ is indeed three-term progression free for any $a\in [0,1]^n$. For  uniformly random $a\in [0,1]^n$, the expected size of this pre-image is $\mathbb{E}[|\varphi_{a,b}^{-1}(S)|]=N\cdot \mu(S)\ge N\cdot 10^{-5}\delta^2n^{-3}\cdot (7/24)^{n/2}$, since for any $x\in \{1,\dots,N\}$ the vector $\varphi_{a,b}(x)$ is distributed uniformly over $[0,1]^n$. Thus, for some choice of $a\in [0,1]^n$ we must have $|\varphi_{a,b}^{-1}(S)|\ge N\cdot 10^{-5}\delta^2n^{-3}\cdot (7/24)^{n/2}$, and so we obtain
\[r_3(N)\ge N\cdot 10^{-5}\delta^2n^{-3}\cdot (7/24)^{n/2}\ge \frac{1}{10^{7}n^3}\cdot N\cdot (7/24)^{n/2}\cdot N^{-2/n}\]
for any choice of $n$. Optimizing $n$, we can again take $n$ as in (\ref{eq-choice-of-n}) and obtain
\begin{align*}
r_3(N)&\ge \frac{1}{10^{7}(2\sqrt{\log_2 N})^{3}}\cdot N\cdot 2^{-\sqrt{\log_2(24/7)}\cdot \sqrt{\log_2 N}-2}\cdot 2^{-\sqrt{\log_2(24/7)}\cdot \sqrt{\log_2 N}}\\
&\ge \frac{N\cdot 2^{-2\sqrt{\log_2(24/7)}\cdot \sqrt{\log_2 N}}}{10^{9}(\log_2 N)^{3/2}}
\end{align*}
if $N$ is sufficiently large. This again shows $r_3(N)\ge N\cdot 2^{-(2\sqrt{\log_2(24/7)}+o(1)) \cdot \sqrt{\log_2 N}}$, with a slightly better bound for the $o(1)$-term compared to the previous proof.

\section{Construction of the building blocks}
\label{sect-constr-building-block}

In this section we prove Proposition~\ref{prop-building-block}, establishing the existence of our desired two-dimensional building blocks for our overall construction. Throughout this section we fix $0<\eps<1$, as in the statement of Proposition~\ref{prop-building-block}.

We now define the desired subset  $T\su [0,1)^2$ as follows.

\begin{definition}\label{defi-T}
Let us define
\begin{align*}
T_1&=\Bigg\lbrace (a,b)\in \Big[\frac{1}{2},1\Big)\times \Big[0,1\Big)\,\Bigg|\,\frac{2}{3}< a+b\le\frac{7}{6}\Bigg\rbrace\\
T_2&=\Bigg\lbrace (a,b)\in \Big[\frac{1}{2},1\Big)\times \Big[0,\frac{1}{2}\Big)\,\Bigg|\,\frac{7}{6}+
\eps\le a+b\le \frac{17}{12}\Bigg\rbrace\\
T_3&=\Bigg\lbrace (a,b)\in \Big[0,\frac{1}{2}\Big)\times \Big[\frac{1}{2},1\Big)\,\Bigg|\,\frac{7}{6}+
\eps\le a+b\le \frac{17}{12} \text{ and }2a+b\ge \frac{3}{2}+\eps\Bigg\rbrace
\end{align*}
and $T=T_1\cup T_2\cup T_3\su [0,1)^2$.
\end{definition}

\begin{figure}
    \centering
\begin{tikzpicture}[scale=6]
\draw[gray!60] (0,0) -- (0,1) -- (1,1) -- (1,0) -- (0,0);
\draw[gray!60] (0,1/2) -- (1,1/2);
\draw[gray!60] (1/2,0) -- (1/2,1);

\def\epsi{0.03}

\draw[gray!20] (0,2/3) -- (2/3,0);
\draw[gray!20] (1/6,1) -- (1,1/6);
\draw[gray!20] (1/6+\epsi,1) -- (1,1/6+\epsi);
\draw[gray!20] (1,5/12) -- (5/12,1);

\draw[->] (0,0) -- (0,1.1);
\draw[->] (0,0) -- (1.1,0);
\draw (0,0) node[anchor=north]{0};
\draw (2/3,0) node[anchor=north]{$\frac{2}{3}$};
\draw (1/2,0) node[anchor=north]{$\frac{1}{2}$};
\draw (1,0) node[anchor=north]{1};
\draw (0,0) node[anchor=east]{0};
\draw (0,1/2) node[anchor=east]{$\frac{1}{2}$};
\draw (0,1) node[anchor=east]{1};
\draw (1,1/6-0.02) node[anchor=west]{\footnotesize{$(1,\frac{1}{6})$}};
\draw (1,1/6+\epsi+0.02) node[anchor=west]{\footnotesize{$(1,\frac{1}{6}+\eps)$}};
\draw (1,5/12) node[anchor=west]{\footnotesize{$(1,\frac{5}{12})$}};
\draw (0.95,1/2) node[anchor=south]{\footnotesize{$(\frac{11}{12},\frac{1}{2})$}};
\draw (0.75,1/2) node[anchor=south]{\footnotesize{$(\frac{2}{3}+\eps,\frac{1}{2})$}};
\draw (1/2,1/6) node[anchor=east]{\footnotesize{$(\frac{1}{2},\frac{1}{6})$}};
\draw (1/2,2/3+0.02) node[anchor=north east]{\footnotesize{$(\frac{1}{2},\frac{2}{3})$}};
\draw (1/2,2/3+\epsi-0.01) node[anchor=south west]{\footnotesize{$(\frac{1}{2},\frac{2}{3}+\eps)$}};
\draw (1/2,11/12) node[anchor=west]{\footnotesize{$(\frac{1}{2},\frac{11}{12})$}};
\draw (5/12,1) node[anchor=south]{\footnotesize{$(\frac{5}{12},1)$}};
\draw (0.2,1) node[anchor=south]{\footnotesize{$(\frac{1}{4}+\frac{\eps}{2},1)$}};
\draw (1/3,5/6+\epsi) node[anchor=east]{\footnotesize{$(\frac{1}{3},\frac{5}{6}+\eps)$}};
\fill[blue!30]  (1,0) -- (1,1/6) -- (1/2,2/3) -- (1/2,1/6) -- (2/3,0) -- (1,0);
\draw[blue,dotted] (1,0) -- (1,1/6);
\draw[blue]  (1,1/6) -- (1/2,2/3) -- (1/2,1/6);
\draw[blue,dotted] (1/2,1/6) -- (2/3,0);
\draw[blue] (2/3,0) -- (1,0);
\fill[blue!30]  (1,1/6+\epsi) -- (1,5/12) -- (11/12,1/2) -- (2/3+\epsi,1/2) -- (1,1/6+\epsi);
\draw[blue,dotted] (1,1/6+\epsi) -- (1,5/12);
\draw[blue] (1,5/12) -- (11/12,1/2);
\draw[blue,dotted] (11/12,1/2) -- (2/3+\epsi,1/2);
\draw[blue] (2/3+\epsi,1/2) -- (1,1/6+\epsi);
\fill[blue!30]  (1/2,11/12) -- (5/12,1) -- (1/4+0.5*\epsi,1) -- (1/3,5/6+\epsi) -- (1/2,2/3+\epsi) --  (1/2,11/12);
\draw[blue] (1/2,11/12) -- (5/12,1);
\draw[blue,dotted] (5/12,1) -- (1/4+0.5*\epsi,1);
\draw[blue] (1/4+0.5*\epsi,1) -- (1/3,5/6+\epsi) -- (1/2,2/3+\epsi);
\draw[blue,dotted] (1/2,2/3+\epsi) --  (1/2,11/12);
\draw[blue] (17/24,5/24) node {$T_1$};
\draw[blue] (9.1/10,4.1/10) node {$T_2$};
\draw[blue] (4.1/10,9.1/10) node {$T_3$};
\end{tikzpicture}
    \caption{The set $T$ defined in Definition \ref{defi-T}}
    \label{figure-set-T}
\end{figure}
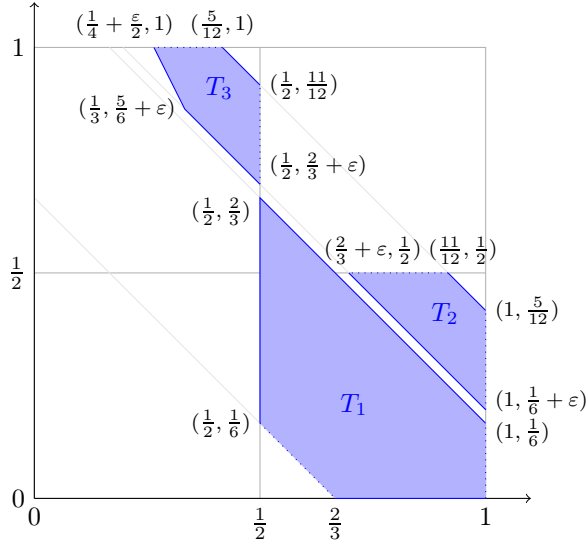

For an illustration of these sets, see Figure \ref{figure-set-T}. Note that the area $\mu(T)$ of the set $T$ can be computed as follows. First, the area of $T_1$ is given by
\[\mu(T_1)=\Big(\frac{1}{2}\Big)^2-\frac{1}{2}\cdot \Big(\frac{1}{3}\Big)^2-\frac{1}{2}\cdot \Big(\frac{1}{6}\Big)^2+\frac{1}{2}\cdot \Big(\frac{1}{6}\Big)^2=\Big(\frac{1}{2}\Big)^2-\frac{1}{2}\cdot \Big(\frac{1}{3}\Big)^2=\frac{9-2}{36}=\frac{7}{36}.\]
Next, the area of $T_2$ is given by
\begin{align*}
\mu(T_2)&=\frac{1}{2}\cdot \Big(\frac{1}{3}\Big)^2-\frac{1}{2}\cdot \Big(\frac{1}{12}\Big)^2-\mu\Big(\{(a,b)\in [1/2,1)\times [0,1/2)\mid 7/6\le a+b< 7/6+\eps\}\Big)\\
&\ge\frac{1}{2}\cdot \Big(\frac{1}{3}\Big)^2-\frac{1}{2}\cdot \Big(\frac{1}{12}\Big)^2-\frac{\eps}{2}= \frac{16-1}{288}-\frac{\eps}{2}= \frac{15}{288}-\frac{\eps}{2}.
\end{align*}
Finally, the area of $T_3$ is given by
\begin{align*}
\mu(T_3)&=\mu\Big(\{(a,b)\in [1/4,1/3)\times [1/2,1)\mid 2a+b\ge 3/2+\eps\}\Big)\\
&\quad\quad\quad+\mu\Big(\{(a,b)\in [1/3,1/2)\times [1/2,1)\mid 7/6+\eps\le a+b\le 17/12\}\Big)\\
&= \frac{1}{2}\cdot \Big(\frac{1}{12}\Big)\cdot \Big(\frac{1}{6}\Big)+\Big(\frac{1}{6}\Big)^2-\frac{1}{2}\cdot \Big(\frac{1}{12}\Big)^2+\frac{1}{2}\cdot \Big(\frac{1}{6}\Big)^2\\
&\quad\quad\quad-\mu\Big(\{(a,b)\in [1/4,1/3)\times [1/2,1)\mid 3/2\le 2a+b< 3/2+\eps\}\Big)\\
&\quad\quad\quad-\mu\Big(\{(a,b)\in [1/3,1/2)\times [1/2,1)\mid 7/6\le a+b\le 7/6+\eps\}\Big)\\
&\ge \frac{2+8-1+4}{288}-\frac{\eps}{12}-\frac{\eps}{6}= \frac{13}{288}-\frac{\eps}{4}.
\end{align*}
Thus, for the area of $T=T_1\cup T_2\cup T_3$ we obtain
\begin{equation}
\mu(T)=\mu(T_1)+\mu(T_2)+\mu(T_3)\ge \frac{7}{36}+\frac{15}{288}-\frac{\eps}{2}+\frac{13}{288}-\frac{\eps}{4}\ge \frac{56+15+13}{288}-\eps=\frac{84}{288}-\eps= \frac{7}{24}-\eps.
\end{equation}

It remains to show that there is a measurable function $f:T\to [0,100/\eps^2]$ such that the condition in Proposition \ref{prop-building-block} is satisfied. To define this function, we first need some more notation. Let us define the function $g:[0,1)\to [0,1/2)$ by setting
\[g(t)=\begin{cases}
t^2 &\text{if }t\in [0,1/2)\\
(t-1/2)^2 &\text{if }t\in [1/2,1)
\end{cases}\]
for all $t\in [0,1)$. Note that this function is measurable and satisfies $0\le g(t)\le (1/2)^2=1/4$ for all $t\in [0,1)$. Furthermore, for any $t,t'\in [0,1)$ with $2t\equiv 2t'\mod 1$, we have $g(t)=g(t')$.

Now, we can define the function $f:T\to [0,100/\eps^2]$ by
\begin{equation}\label{eq-defi-f}
f(x)=\frac{24}{\eps^2}\cdot (x_1+x_2)^2+6\cdot g(x_1)
\end{equation}
for all $x\in T\su [0,1)^2$, where we write the coordinates of $x$ as $x=(x_1,x_2)$ as usual. Note that $f$ is a measurable function and we indeed have $0\le f(x)\le (24/\eps^2) \cdot 2^2+6\cdot 1/4\le 100/\eps^2$ for all $x\in T$. We need to show that for any $x,y,z\in T$ with $x+z\equiv 2y\mod 1$ we have
\begin{equation}\label{eq-inequality-to-show}
f(x) + f(z) \ge 2f(y) + (x_1-z_1)^2 + (x_2-z_2)^2.
\end{equation}

In order to show this, we begin with some simple facts.

\begin{fact}\label{fact-simple-inequalities}
For any $s,t\in \mathbb{R}$, the following statements hold:
\begin{itemize}
\item[(i)] $s^2+t^2=2\cdot (s/2+t/2)^2+(s-t)^2/2$.
\item[(ii)] $(s+t)^2+t^2\ge s^2/2$.
\item[(iii)] If $t\ge 0$, then $(t+1/2)^2\ge t^2+1/4$.
\end{itemize} 
\end{fact}
\begin{proof}
For (i), we observe
\[s^2+t^2=2\cdot \Big(\frac{s+t}{2}\Big)^2+2\cdot \Big(\frac{s-t}{2}\Big)^2=2\cdot \Big(\frac{s}{2}+\frac{t}{2}\Big)^2+ \frac{(s-t)^2}{2}.\]
Plugging $s+t$ and $t$ into this equation, we obtain
\[(s+t)^2+t^2=2\cdot \Big(\frac{s+t}{2}+\frac{t}{2}\Big)^2+ \frac{s^2}{2}\ge \frac{s^2}{2},\]
showing (ii). For (iii), we simply note that
\[\Big(t+\frac{1}{2}\Big)^2=t^2+t+\frac{1}{4}\ge t^2+\frac{1}{4}\]
if $t\ge 0$.
\end{proof}

\begin{fact}\label{fact-simple-properties-coordinate-sum}
For the set $T\su [0,1)^2$ in Definition \ref{defi-T}, the following statements hold:
\begin{itemize}
\item[(i)] For all $(a,b)\in T$, we have $2/3<a+b\le 17/12$.
\item[(ii)] For all $(a,b)\in T$, we have $g(a)\ge (a-1/2)^2$.
\item[(iii)] For all $(a,b), (a',b') \in T$ with $a+a'< 1$, we have $a+b+a'+b'> 11/6$.
\end{itemize} 
\end{fact}

\begin{proof}
Note that (i) follows immediately from the definition of $T$.

To show (ii), first note that $g(a)=(a-1/2)^2$ for all $a\in [1/2,1)$. So we have $g(a)= (a-1/2)^2$ for all $(a,b)\in T_1\cup T_2$. It remains to consider $(a,b)\in T_3$, then $g(a)=a^2$. Furthermore, by $2a+b\ge 3/2+\eps>3/2$ and $b<1$, we have $a> 1/4$ and hence $g(a)=a^2> a^2-a+1/4=(a-1/2)^2$.

It remains to show (iii), so consider $(a,b), (a',b') \in T$ with $a+a'< 1$. Then we must have $a< 1/2$ or $a'<1/2$, so let us assume $a'<1/2$ without loss of generality. This implies $(a',b')\in T_3$ and hence $a'+b'\ge 7/6+\eps> 7/6$ by the definition of $T_3$. Combining this with the lower bound in (i), we obtain
\[a+b+a'+b'> \frac{2}{3}+\frac{7}{6}=\frac{11}{6}.\qedhere\]
\end{proof}

We can now derive various consequences for points $ x,y,z\in T$ such that $x+z\equiv 2y\mod 1$, with the aim of proving (\ref{eq-inequality-to-show}).

\begin{lemma}\label{lemma-AP-coordinate-sum-simple}
Let $T\su [0,1)^2$ be defined as in Definition \ref{defi-T}. Consider points $ x,y,z\in T$ such that $x+z\equiv 2y\mod 1$. Then we have either
\[y_1+y_2= \frac{x_1+x_2}{2}+\frac{z_1+z_2}{2}\]
or 
\[y_1+y_2= \frac{x_1+x_2}{2}+\frac{z_1+z_2}{2}-\frac{1}{2}\]
\end{lemma}
\begin{proof}
Let $y'=(y_1',y_2')\in [0,1)^2$  be the mid-point of $x=(x_1,x_2)$ and $z=(z_1,z_2)$ (given by $y_1'=(x_1+z_1)/2$ and $y_2'=(x_2+z_2)/2$). Note that
\[(2y_1,2y_2)=2y\equiv x+z=(x_1,x_2)+(z_1,z_2)=(x_1+z_1,x_2+z_2)=(2y_1',2y_2')\mod 1,\]
meaning that $2y_1\equiv 2y_1'\mod 1$ and $2y_2\equiv 2y_2'\mod 1$. Hence $y_1-y_1'\in \{-1/2,0,1/2\}$ and $y_2-y_2'\in \{-1/2,0,1/2\}$. Since $(y_1,y_2)\not\in [0,1/2)^2$ (as $T\cap [0,1/2)^2=\emptyset$), we cannot have both $y_1-y_1'=-1/2$ and  $y_2-y_2'=-1/2$. So we can conclude that $(y_1+y_2)-(y_1'+y_2')\in \{-1/2,0,1/2,1\}$.

The assertion of the lemma is that $(y_1+y_2)-(y_1'+y_2')\in \{-1/2,0\}$. Hence it suffices to prove $(y_1+y_2)-(y_1'+y_2')<1/2$. So let us assume for contradiction that $(y_1+y_2)-(y_1'+y_2')\ge 1/2$.

Since $y_2-y_2'\le 1/2$, we must therefore have $y_1\ge y_1'$. Similarly, we must have $y_2\ge y_2'$, since $y_1-y_1'\le 1/2$.

If $y_1'<1/2$ (i.e. if $x_1+z_1<1$), by Fact \ref{fact-simple-properties-coordinate-sum}(iii) we have
\[y_1'+y_2'=\frac{(x_1+x_2)+(z_1+z_2)}{2}> \frac{11/6}{2}=\frac{11}{12}=\frac{17}{12}-\frac{1}{2}\ge y_1+y_2-\frac{1}{2},\]
where the last inequality follows from Fact \ref{fact-simple-properties-coordinate-sum}(i). This contradicts $(y_1+y_2)-(y_1'+y_2')\ge 1/2$, so we must have $y_1'\ge 1/2$.

Now, we obtain $1/2\le y_1'\le y_1<1$. Together with $y_1-y_1'\in \{-1/2,0,1/2\}$, this implies $y_1=y_1'$. So from $(y_1+y_2)-(y_1'+y_2')\ge 1/2$, we obtain $y_2\ge y_2'+1/2\ge 1/2$.

Thus, we have $y_1=y_1'\ge 1/2$ and $y_2\ge 1/2$, meaning that $(y_1,y_2)\in [1/2,1)^2$. By the definition of $T$, this implies $y_1+y_2\le 7/6$ and hence
\[y_1'+y_2'=\frac{(x_1+x_2)+(z_1+z_2)}{2}> \frac{2/3+2/3}{2}=\frac{2}{3}=\frac{7}{6}-\frac{1}{2}\ge y_1+y_2-\frac{1}{2},\]
where the first inequality follows from Fact \ref{fact-simple-properties-coordinate-sum}(i). This is again a contradiction to $(y_1+y_2)-(y_1'+y_2')\ge 1/2$. 
\end{proof}

\begin{corollary}\label{coro-AP-coordinate-sum}
Let $T\su [0,1)^2$ be defined as in Definition \ref{defi-T}, and let $ x,y,z\in T$ be points such that $x+z\equiv 2y\mod 1$. Then at least one of the following two statements holds
\begin{itemize}
\item[(a)] $(x_1+x_2)^2+(z_1+z_2)^2\ge 2\cdot (y_1+y_2)^2+\eps^2/2$.
\item[(b)] $|x_1+x_2-z_1-z_2|<\eps$ and $(x_1+x_2)^2+(z_1+z_2)^2= 2\cdot (y_1+y_2)^2+(x_1+x_2-z_1-z_2)^2/2$.
\end{itemize} 
\end{corollary}

\begin{proof}
By Lemma \ref{lemma-AP-coordinate-sum-simple}, either $y_1+y_2=(x_1+x_2)/2+(z_1+z_2)/2$ or $y_1+y_2=(x_1+x_2)/2+(z_1+z_2)/2-1/2$ must hold. In the latter case,  we have
\[(x_1+x_2)^2+(z_1+z_2)^2\ge 2\cdot \Big(\frac{x_1+x_2}{2}+\frac{z_1+z_2}{2}\Big)^2=2\cdot (y_1+y_2+1/2)^2\ge 2\cdot (y_1+y_2)^2+1/2\ge 2\cdot (y_1+y_2)^2+\eps^2/2,\]
where the first inequality follows from Fact \ref{fact-simple-inequalities}(i) and the second inequality from  Fact \ref{fact-simple-inequalities}(iii). So let us now assume that $y_1+y_2=(x_1+x_2)/2+(z_1+z_2)/2$, then by Fact \ref{fact-simple-inequalities}(i) we have
\[(x_1+x_2)^2+(z_1+z_2)^2= 2\cdot \Big(\frac{x_1+x_2}{2}+\frac{z_1+z_2}{2}\Big)^2+(x_1+x_2-z_1-z_2)^2/2=2\cdot (y_1+y_2)^2+(x_1+x_2-z_1-z_2)^2/2.\]
If $|x_1+x_2-z_1-z_2|\ge \eps$, this implies the inequality in (a). And if $|x_1+x_2-z_1-z_2|<\eps$, then (b) holds.
\end{proof}

Recall that our goal is to prove (\ref{eq-inequality-to-show}) for any $ x,y,z\in T$ with $x+z\equiv  2y\mod 1$. If (a) in Corollary \ref{coro-AP-coordinate-sum} holds, then it is not hard to show inequality (\ref{eq-inequality-to-show}). The next lemma helps with showing the inequality in case (b) in Corollary \ref{coro-AP-coordinate-sum} holds.

\begin{lemma}\label{lem-second-property-set-T}
Let $T\su [0,1)^2$ be defined as in Definition \ref{defi-T}, and let $x,z\in T$ be such that $|x_1+x_2-z_1-z_2|<\eps$. Suppose that $x_1\ge 1/2$ or $z_1\ge 1/2$. Then we have $x_1+z_1\ge 1$.
\end{lemma}

\begin{proof}
Recall that for all points $(a,b)\in T_1$ we have $a+b\le 7/6$, whereas for all points $(a,b)\in T_2\cup T_3$ we have $a+b\ge 7/6+\eps$. Since $|x_1+x_2-z_1-z_2|<\eps$, this means that we either have $x,z\in T_1$ or $x,z\in T_2\cup T_3$.

If $x,z\in T_1$, then we have $x_1\ge 1/2$ and $z_1\ge 1/2$ and hence $x_1+z_1\ge 1$, as desired. So let us assume that $x,z\in T_2\cup T_3$. If $x,z\in T_2$, we similarly have  $x_1\ge 1/2$ and $z_1\ge 1/2$ and hence $x_1+z_1\ge 1$. If $x,z\in T_3$, then we would have $x_1< 1/2$ and $z_1< 1/2$, which contradicts the assumption in the lemma.

So it only remains to consider the case that one of the points $x$ and $z$ is contained in the set $T_2$, and the other point in the set $T_3$. Let us assume without loss of generality that $x\in T_2$ and $z\in T_3$. Then we have
\[x_1+z_1=x_1+2z_1+z_2-(z_1+z_2)\ge x_1+3/2+\eps-(z_1+z_2)\ge x_1+3/2+\eps-(x_1+x_2+\eps)=3/2-x_2\ge 1,\]
where the first inequality follows from the definition of $T_3$, the second inequality follows from the assumption $|x_1+x_2-z_1-z_2|<\eps$, and the third inequality follows from the definition of $T_2$.
\end{proof}

\begin{corollary}\label{coro-option-b}
Let $T\su [0,1)^2$ be defined as in Definition \ref{defi-T}, and let $x,y,z\in T$ be such that $x+z\equiv 2y\mod 1$ and $|x_1+x_2-z_1-z_2|<\eps$. Then we have
\[g(x_1)+g(z_1)\ge 2\cdot g(y_1)+(x_1-z_1)^2/2.\]
\end{corollary}
\begin{proof}
Since $x+z\equiv 2y\mod 1$, we have $x_1+z_1\equiv 2y_1\mod 1$ and hence $2(x_1/2+z_1/2)\equiv 2y_1\mod 1$. So we have $g(y_1)=g(x_1/2+z_1/2)$, and therefore the claimed inequality is equivalent to
\[g(x_1)+g(z_1)\ge 2\cdot g(x_1/2+z_1/2)+(x_1-z_1)^2/2.\]

If $x_1\ge 1/2$ or $z_1\ge 1/2$, then by Lemma \ref{lem-second-property-set-T} we have $x_1/2+z_1/2\ge 1/2$ and hence $g(x_1/2+z_1/2)=(x_1/2+z_1/2-1/2)^2$. Thus, by Fact \ref{fact-simple-properties-coordinate-sum}(ii) and Fact \ref{fact-simple-inequalities}(i), in this case we obtain
\[g(x_1)+g(z_1)\ge (x_1-1/2)^2+(z_1-1/2)^2=2\cdot (x_1/2+z_1/2-1/2)^2+(x_1-z_1)^2/2=2\cdot g(x_1/2+z_1/2)+(x_1-z_1)^2/2,\]
as desired. 

It remains to consider the case that $x_1< 1/2$ and $z_1< 1/2$, then we also have $x_1/2+z_1/2<1/2$. Thus, observing that $g(x_1)=x_1^2$ and $g(z_1)=z_1^2$ and $g(x_1/2+z_1/2)=(x_1/2+z_1/2)^2$, by Fact \ref{fact-simple-inequalities}(i) we obtain
\[g(x_1)+g(z_1)= x_1^2+z_1^2=2\cdot (x_1/2+z_1/2)^2+(x_1-z_1)^2/2=2\cdot g(x_1/2+z_1/2)+(x_1-z_1)^2/2.\qedhere\]
\end{proof}

Finally, we are ready to prove Proposition \ref{prop-building-block}.

\begin{proof}[Proof of Proposition \ref{prop-building-block}]
Let $T\su [0,1)^2$ be given as in Definition \ref{defi-T}, and let $f:T\to [0,100/\eps^2]$ be the measurable function defined in (\ref{eq-defi-f}). In order to prove the proposition it remains to show the inequality
\[f(x) + f(z) \ge 2\cdot f(y) + (x_1-z_1)^2 + (x_2-z_2)^2\]
for $x,y,z\in T\su [0,1)^2$ with $x+z\equiv 2y \mod 1$. For any such $x,y,z\in T$, one of the two statements in Corollary \ref{coro-AP-coordinate-sum} holds. If (a) holds, then we have
\begin{align*}
f(x)+f(z)&=\frac{24}{\eps^2}\cdot (x_1+x_2)^2+6\cdot g(x_1)+\frac{24}{\eps^2}\cdot (z_1+z_2)^2+6\cdot g(z_1)\\
&\ge\frac{24}{\eps^2}\cdot\Big((x_1+x_2)^2+(z_1+z_2)^2\Big)\\
&\ge\frac{24}{\eps^2}\cdot\Big(2\cdot (y_1+y_2)^2+\eps^2/2\Big)\\
&=2\cdot \frac{24}{\eps^2}\cdot(y_1+y_2)^2+12\\
&\ge 2\cdot \frac{24}{\eps^2}\cdot(y_1+y_2)^2+2\cdot 6\cdot (1/4)+1+1\\
&\ge 2\cdot \frac{24}{\eps^2}\cdot(y_1+y_2)^2+2\cdot 6\cdot g(y_1)+(x_1-z_1)^2+(x_2-z_2)^2\\
&= 2\cdot f(y) + (x_1-z_1)^2 + (x_2-z_2)^2.
\end{align*}
If (b) in Corollary \ref{coro-AP-coordinate-sum} holds, then in particular by Corollary \ref{coro-option-b} we have
\[g(x_1)+g(z_1)\ge 2\cdot g(y_1)+(x_1-z_1)^2/2.\]
Together with the equation in (b), this yields
\begin{align*}
f(x)+f(z)&=\frac{24}{\eps^2}\cdot (x_1+x_2)^2+6\cdot g(x_1)+\frac{24}{\eps^2}\cdot (z_1+z_2)^2+6\cdot g(z_1)\\
&=\frac{24}{\eps^2}\cdot\Big((x_1+x_2)^2+(z_1+z_2)^2\Big)+6\cdot \Big(g(x_1)+g(z_1)\Big)\\
&\ge\frac{24}{\eps^2}\cdot\Big(2\cdot (y_1+y_2)^2+(x_1+x_2-z_1-z_2)^2/2\Big)+6\cdot \Big(2\cdot g(y_1)+(x_1-z_1)^2/2\Big)\\
&=2\cdot \Big(\frac{24}{\eps^2}\cdot(y_1+y_2)^2+6\cdot g(y_1)\Big)+\frac{12}{\eps^2}\cdot (x_1+x_2-z_1-z_2)^2+3\cdot (x_1-z_1)^2\\
&\ge 2\cdot f(y) + (x_1-z_1)^2 + 2\cdot (x_1+x_2-z_1-z_2)^2+2\cdot (x_1-z_1)^2\\
&\ge 2\cdot f(y) + (x_1-z_1)^2 + (x_2-z_2)^2,
\end{align*}
where in the last step we used Fact \ref{fact-simple-inequalities}(ii). In either case we proved the desired inequality, finishing the proof of Proposition \ref{prop-building-block}.
\end{proof}

\end{document}